\documentclass [12pt,a4paper,reqno]{amsart}
\input amssymb.sty
\def\Xint#1{\mathchoice
   {\XXint\displaystyle\textstyle{#1}}%
   {\XXint\textstyle\scriptstyle{#1}}%
   {\XXint\scriptstyle\scriptscriptstyle{#1}}%
   {\XXint\scriptscriptstyle\scriptscriptstyle{#1}}%
   \!\int}
\def\XXint#1#2#3{{\setbox0=\hbox{$#1{#2#3}{\int}$}
     \vcenter{\hbox{$#2#3$}}\kern-.5\wd0}}

\def\dashint{\Xint-}

\usepackage{amsmath,amsthm}
\usepackage{amsfonts}
\allowdisplaybreaks
\renewcommand{\a} {\alpha}
\renewcommand{\b} {\beta}

\renewcommand{\d} {\delta}

\newcommand{\g} {\gamma}
\renewcommand{\l} {\lambda}
\renewcommand{\S} {\Sigma}

\newcommand {\dist} {{\rm dist}}

\newcounter{lemma}[section]

\newcounter{corol}[section]

\newcounter{rem}[section]

\newcounter{theo}[section]

\newcounter{propo}[section]

\begin{document}
\title{Topological Mappings with Controlled $p$-Moduli}
\author{Anatoly Golberg and Ruslan Salimov}
\begin{abstract}
We study homeomorphisms of controlled $p$-module by certain
integrals. In this way, we establish various properties of mappings
and show that their features are close to quasiconformal and
bilipschitz mappings.

\end{abstract}

\date{\today \hskip 4mm (\texttt{TMCM.tex})}
\maketitle


\bigskip
{\small {\textbf {Key words:} $p$-module of $k$-dimensional surface
families, weighted $p$-module, ring $Q$-homeomor\-phisms, lower
$Q$-homeomorphisms, $\alpha$-inner and $\alpha$-outer dilatations,
bilipschitz mappings

\bigskip

\medskip

\section{Introduction}
A great interest to studying various
classes of homeomorphisms and more general mappings is motivated by
needs of many fields in Modern Mathematics. Some of basic classes
are close to quasiconformal and bilipschitz homeomorphisms. A main
characterization of such mappings is obtained by extension of
quasiinvariance of the conformal moduli and $p$-moduli via
inequalities containing integrals depending on a given measurable
functions and admissible metrics (cf. \cite{Gol09, Gol11, KR08,
MRSY09, Sal08, Sal11, SS10}). Such representation of moduli can be
treated as the quasiinvariance of weighted moduli (cf. \cite{AC71,
Cris11, Tam98}).

Let $f:\,G\to G^*$, $G, G^*\subset\mathbb R^n$, be a homeomorphism
such that $f$ and $f^{-1}$ are differentiable almost everywhere
(a.e.) with nonzero Jacobians in $G$ and $G^*$, respectively. It was
shown in \cite{Gol11}, that under more restrictive conditions on $f$
the following bounds for the $\alpha$-module of $k$-dimensional
surface families
\begin{equation}\label{eq1.1}
\inf\limits_{\varrho\in{\rm extadm}\;\mathcal
S_k}\int\limits_G\frac{\varrho^\alpha(x)} {H_{O,\alpha}(x,f)}\,dx\le
\mathcal M_\alpha(f(\mathcal S_k))\le\inf\limits_{\rho\in{\rm
adm}\;\mathcal S_k}\int\limits_G\rho^\alpha(x) H_{I,\alpha}(x,f)\,dx
\end{equation}
are fulfilled. Here $H_{I,\alpha}(x,f)$ and $H_{O,\alpha}(x,f)$
stand for the $\alpha$-inner and $\alpha$-outer dilatations of $f$
at $x\in G$ (see, e.g., \cite{Gol05}).

In the case, when these dilatations are bounded, i.e.
$H_{I,\alpha}(x,f)\le K$ and $H_{O,\alpha}(x,f)\le K$ with some
absolute constant $K$ in $G$, one obtains the well-known class of
bilipschitz mappings in $G$ (cf. \cite{Geh71, Gol05, GR90}).

In this paper, we consider the homeomorphisms satisfying at least
one of the following conditions
\begin{equation}\label{eq1.2}
\mathcal M_\alpha(f(\mathcal S_k))\le\inf\limits_{\rho\in{\rm
adm}\;\mathcal S_k}\int\limits_G\rho^\alpha(x) Q(x)\,dx,
\end{equation}
\begin{equation}\label{eq1.3}
\mathcal M_\alpha(f(\mathcal S_k))\ge\inf\limits_{\varrho\in{\rm
extadm}\;\mathcal S_k}\int\limits_G\frac{\varrho^\alpha(x)}
{Q(x)}\,dx,
\end{equation}
with a given measurable function $Q:G\to [0,\infty]$. For such
mappings the problem can be formulated somethat similarly to the
classical problem on the properties of solutions to the Beltrami
equation $f_{\bar z}=\mu(z)f_z$, for which the properties of $f$ are
investigated in their dependence on the features of $\mu$.

The main cases in (\ref{eq1.2})-(\ref{eq1.3}) relate to $k=1$ and
$k=n-1$, i.e. to moduli of curve and of $(n-1)$-surface families. We
show that inequality (\ref{eq1.2}) yields differentiability a.e.,
the $(N)$-property, boundedness of the $\alpha$-inner dilatation. We
also provide the necessary and sufficient condition for a
homeomorphism to satisfy (\ref{eq1.3}). Finally, we establish the
relationship between homeomorphisms satisfying (\ref{eq1.2}) for
$k=1$ and (\ref{eq1.3}) for $k=n-1$.

\section{Dilatations in $\mathbb R^n$}
Let $A:\;\mathbb{R}^n\to
\mathbb{R}^n$ be a linear bijection. The numbers

\begin{equation*}
H_{I,\alpha}(A)=\frac{|\det A|}{l^\alpha(A)},\quad
H_{O,\alpha}(A)=\frac{||A||^\alpha}{|\det A|}, \quad \alpha\ge 1,
\end{equation*}
are called the \textit{$\alpha$-inner} and \textit{$\alpha$-outer}
dilatations of $A$, respectively. Here
\begin{equation*}
l(A)=\min\limits_{|h|=1} |Ah|,\ \ ||A||=\max\limits_{|h|=1} |Ah|,
\end{equation*}
denote the minimal and maximal stretching of $A$ and $\det A$ is the
determinant of $A$.

Let $G$ and $G^*$ be two bounded domains in $\mathbb{R}^n$, $n \geq
2$, and let a mapping $f:\;G\to G^*$ be differentiable at a point $x
\in G$. This means there exists a linear mapping
$f^{\prime}(x):\;\mathbb{R}^n\to \mathbb{R}^n$, called the (strong)
derivative of the mapping $f$ at $x$, such that
\begin{equation*}
f(x+h)=f(x)+f^{\prime}(x)h+ \omega (x,h) |h|,
\end{equation*}
where $\omega (x,h) \to 0$ as $h \to 0$.

We denote
\begin{equation*}
H_{I,\alpha}(x,f)=H_{I,\alpha}(f^{\prime}(x)),\quad
H_{O,\alpha}(x,f)=H_{O,\alpha}(f^{\prime}(x)).
\end{equation*}
These quantities naturally extend the classical quasiconformal
dilatations (inner and outer) by
\begin{equation*}
H_I(x,f)=H_{I,n}(x,f),\quad H_O(x,f)=H_{O,n}(x,f).
\end{equation*}
The third dilatation of quasiconformality called \textit{linear}
\begin{equation*}
H(x,f)=\frac{||f^{\prime}(x)||}{l(f^{\prime}(x))}
\end{equation*}
is a direct analog of the classical planar Lavrentiev dilatation.
The $\alpha$-inner and $\alpha$-outer dilatations provide a class of
mappings whose basic properties are close to quasiconformal
homeomorphisms. On the other hand, there are some essential
differences caused by the fact that the dilatations $H_I(x,f)$ and
$H_O(x,f)$ are always greater than or equal to 1, while the
$\alpha$-inner and outer dilatations range between 0 and $\infty$.

\medskip
We consider the homeomorphisms $f$ which are differentiable almost
everywhere in $G$, and fix the real numbers $\alpha, \beta$
satisfying $1\le\alpha<\beta<\infty$. Define
\begin{equation*}
HI_{\alpha,\beta}(f)=\int\limits_G H_{I,\alpha}^
\frac\beta{\beta-\alpha}(x,f)\,dx,\quad
HO_{\alpha,\beta}(f)=\int\limits_G H_{O,\beta}^
\frac\alpha{\beta-\alpha}(x,f)\,dx,
\end{equation*}
and call these quantities the \textit{inner} and \textit{outer mean
dilatations} of a mapping $f:\;G\to\mathbb R^n$ in $G$.

\medskip
Define for the fixed real numbers $\alpha,\beta,\gamma,\delta$ such
that $1\le\alpha<\beta<\infty$, $1\le\gamma<\delta<\infty$, the
\textit{class of mappings with finite mean dilatations} which
consists of homeomorphisms $f:\;G\to G^*$ satisfying:

\noindent (i) $f$ and $f^{-1}$ are in $W^{1,1}_{\rm loc}$,

\noindent (ii) $f$ and $f^{-1}$ are differentiable, with Jacobians
$J(x,f)\ne 0$ and $J(y,f^{-1})\ne 0$ a.e. in $G$ and $G^*$,
respectively,

\noindent (iii) the inner and the outer mean dilatations
$HI_{\alpha,\beta}(f)$ and $HO_{\gamma,\delta}(f)$ are finite.

The mappings with finite mean dilatations were investigated in \cite
{Gol05}.

The relations between the classical quasiconformal dilatations
\begin{equation*}
\begin{split}
H(x,f)&\le \min(H_I(x,f),H_O(x,f))\le H^{n/2}(x,f) \\
&\le\max(H_I(x,f),H_O(x,f))\le H^{n-1}(x,f).
\end{split}
\end{equation*}
show that they are finite or infinity simultaneously. However, this
needs not be true for the mean dilatations.

The following example shows that the unboundedness of one from these
dilatations does not depend on the value of another mean dilatation.

\medskip
Consider the unit cube
\begin{equation*}
G=\{x=(x_1,\ldots,x_n):0<x_k<1,k=1,\ldots,n\}
\end{equation*}
and let
\begin{equation*}
f(x)=\biggl(x_1,\ldots,x_{n-1},\frac{x_n^{1-c}}{1-c}\biggr), \quad
0<c<1.
\end{equation*}
An easy computation shows that $f$ belongs to the class of mappings
with finite mean dilatations if and only if
\begin{equation*}
0<c<1-\alpha/\beta \quad \text{and} \quad 0<c<1-(\g-1)\d/(\d-1)\g.
\end{equation*}
When
\begin{equation*}
1-\a/\b\le c<1 \quad \text{and} \quad 1-(\g-1)\d/(\d-1)\g\le c<1
\end{equation*}
we have $HI_{\alpha,\beta}(f)=\infty$ and
$HO_{\gamma,\delta}(f)=\infty$, respectively. Thus, by suitable
choice of parameters $c,\alpha,\beta,\gamma,\delta$, one obtains any
desired relations between $HI_{\alpha,\beta}(f)$ and
$HO_{\gamma,\delta}(f)$.

\section{$\alpha$-moduli of $k$-dimensional surfaces and related
classes of homeomorphisms}
Now we give a geometric (modular)
description of quasiconformality in $\mathbb R^n$ starting with the
definition of $k$-dimensional Hausdorff measure $H^k$,
$k=1,\ldots,n-1$ in $\mathbb R^n$. For a given $E\subset\mathbb
R^n$, put
\begin{equation*}
H^k(E)=\sup\limits_{r>0} H^k_r(E),
\end{equation*}
where
\begin{equation*}
H^k_r(E)=\Omega_k\inf\sum\limits_i (\delta_i/2)^k.
\end{equation*}
Here the infimum is taken over all countable coverings $\{E_i,
i=1,2,\ldots \}$ of $E$ with diameters $\delta_i$, and $\Omega_k$ is
the volume of the unit ball in $\mathbb R^k$.

Let $\mathcal{S}$ be a $k$-dimensional surface, which means that
$\mathcal{S}:\;D_s\to \mathbb{R}^n$ is a continuous image of the
closed domain $D_s\subset \mathbb{R}^k$. We denote by
\begin{equation*}
N(\mathcal S,y)={\rm card}\, \mathcal S^{-1}(y)={\rm card} \{x\in
D_s: \mathcal S(x)=y\}
\end{equation*}
the multiplicity function of the surface $\mathcal S$ at the point
$y\in\mathbb R^n$. For a given Borel set $B\subseteq \mathbb R^n$,
the $k$-dimensional Hausdorff area of $B$ in $\mathbb R^n$
associated with the surface $\mathcal S$ is determined by
\begin{equation*}
\mathcal H_\mathcal S(B)=\mathcal H_\mathcal S^k(B)= \int\limits_B
N(\mathcal S,y)\, d H^k y.
\end{equation*}
If $\rho:\,\mathbb R^n\to [0,\infty]$ is a Borel function, the
integral of $\rho$ over $\mathcal S$ is defined by
\begin{equation*}
\int\limits_\mathcal S \rho\, d\sigma_k=\int\limits_{\mathbb R^n}
\rho(y) N(\mathcal S,y)\,dH^k y.
\end{equation*}

Let $\mathcal{S}_k$ be a family of $k$-dimensional surfaces
$\mathcal{S}$ in $\mathbb{R}^n$, $1\le k\le n-1$ (curves for $k=1$).
The \textit{$\alpha$-module} of $\mathcal{S}_k$ is defined as
\begin{equation*}
\mathcal M_\alpha(\mathcal{S}_k)=\inf \int\limits_{\mathbb{R}^n}
\rho^\alpha\,dx,\ \ \alpha\ge k,
\end{equation*}
where the infimum is taken over all Borel measurable functions
$\rho\ge 0$ and such that
\begin{equation*}
\int\limits_\mathcal S\rho^k\,d\sigma_k\ge 1
\end{equation*}
for every $\mathcal{S}\in \mathcal{S}_k$. We call each such $\rho$
an \textit{admissible function} for $\mathcal{S}_k$ ($\rho\in {\rm
adm}\, \mathcal S_k$). The $n$-module $\mathcal M_n( \mathcal{S}_k)$
will be denoted by $\mathcal M(\mathcal{S}_k)$.

Following \cite{KR08}, a metric $\rho$ is said to be
\textit{extensively admissible} for $\mathcal{S}_k$ ($\rho\in {\rm
extadm}\, \mathcal S_k$) with respect to $\alpha$-module if $\rho\in
{\rm adm}\, (\mathcal S_k\backslash\widetilde{\mathcal S}_k)$ such
that $\mathcal M_\alpha(\widetilde{\mathcal{S}}_k)=0$.

Accordingly, we say that a property $P$ holds for almost every
$k$-dimensional surface, if $P$ holds for all surfaces except a
family of zero $\alpha$-module.

\medskip
We also remind that a continuous mapping $f$ satisfies
$(N)$-property with respect to $k$-dimensional Hausdorff area if
$\mathcal H_\mathcal S^k(f(B))=0$ whenever $\mathcal H_\mathcal
S^k(B)=0$. Similarly, $f$ has $(N^{-1})$-property if $\mathcal
H_\mathcal S^k(B)=0$ whenever $\mathcal H_\mathcal S^k(f(B))=0$.

\medskip
Now we provide the bounds for the $\alpha$-module of $k$-dimensional
surfaces (see \cite{Gol11}).

\medskip
\begin{theo}
{\em Let $f:\;G\to \mathbb R^n$ be a homeomorphism satisfying (i)-(ii)
with $H_{I,\alpha}, H^{-1}_{O,\alpha}\in L^1_{\rm loc}(G)$. Suppose
that for some $k$, $1\le k\le n-1$ $(k\le\alpha)$, and for almost
every $k$-dimensional surface $\mathcal S$ and its image $\mathcal
S^*$ ($\mathcal S=f^{-1}(\mathcal S^*)$) the restriction
$f|_{\mathcal S}$ has the $(N)$ and $(N^{-1})$-properties with
respect to $k$-dimensional Hausdorff area in $G$ and $G^*=f(G)$,
respectively. Then the double inequality (\ref{eq1.1}) holds for any
family $\mathcal S_k$ of $k$-dimensional surfaces in $G$, and for
each $\rho\in{\rm adm}\,\mathcal S_k$ and $\varrho\in{\rm
extadm}\,\mathcal S_k$ with respect to the $\alpha$-module.}
\end{theo}

\medskip
Further we use the following lemma from \cite{MRSY09}.

\medskip
\begin{lemma}\label{lem}
{\em Let $(X,\mu)$ be s measure space with finite measure $\mu$, and let
$\varphi:G\to (0,\infty)$ be a measurable function. Set
\begin{equation*}
I(\varphi,\alpha)=\inf\limits_\rho\int\limits_X \varphi\rho^\alpha\;
d\mu,\quad 1<\alpha<\infty,
\end{equation*}
where the infimum is taken over all Borel nonnegative measurable
functions $\rho:\,X\to [0,\infty]$ satisfying
$\int\limits_X\rho\,d\mu=1$. Then
\begin{equation*}
I(\varphi,\alpha)=\biggl(\int\limits_X\varphi^\frac{1}{1-\alpha}\;d\mu\biggr)^{1-\alpha}
\end{equation*}
and the infimum is attained only for the metric
\begin{equation*}
\rho=\left(\int\limits_X
\varphi^\frac{1}{1-\alpha}\;d\mu\right)^{-1}\varphi^\frac{1}{1-\alpha}.
\end{equation*}}
\end{lemma}

\medskip

Throughout the paper, we use the following notations. A ring domain
$\mathcal{R}\subset\mathbb{R}^n$ is a bounded domain whose
complement consists of two components $C_0$ and $C_1$. The sets
$F_0=\partial C_0$ and $F_1=\partial C_1$ are two boundary
components of $D$. We assume for definiteness that $\infty\in C_1$.

We say that a curve $\gamma$ \textit{joins the boundary components
in} $\mathcal{R}$ if $\gamma$ is located in $\mathcal{R}$, except
for its endpoints, one of which lies on $F_0$ and the second on
$F_1$. A compact set $\S $ is said to \textit{separate the boundary
components of} $\mathcal{R}$ if $\S \subset \mathcal{R}$ and if
$C_0$ and $C_1$ are located in different components of $C\Sigma$.
Denote by $\Gamma_\mathcal{R}$ the family of all locally rectifiable
curves $\gamma$ which join the boundary components of $\mathcal{R}$
and by $\Sigma_\mathcal{R}$ the family of all compact piecewise
smooth $(n-1)$-dimensional surfaces $\Sigma$ which separate the
boundary components of $\mathcal{R}$.

\medskip
The following relation
\begin{equation}\label{ziem}
\mathcal M_p(\Gamma_\mathcal R)=\frac{1}{\mathcal
M^{p-1}_\alpha(\Sigma_\mathcal R)},\qquad
\alpha=\frac{p(n-1)}{p-1},\qquad 1<p<\infty, \qquad
n-1<\alpha<\infty,
\end{equation}
between the $p$-moduli of $\Sigma_\mathcal R$ and $ \Gamma_\mathcal
R$ follows from the results of Ziemer \cite{Zie69} and Hesse
\cite{Hes75} on the moduli and the extremal lengths. Observe that
$p$-moduli ($p\ne n$) are not conformal invariants even under linear
mappings, i.e. such mappings do not preserve the value of
$p$-module.

The $p$-module of a spherical ring $A(x_0;a,b)=\{x\in\mathbb{R}^n:
0<a<|x-x_0|<b\}$ is equal to
\begin{equation*}
\mathcal
M_p(\Gamma_A)=\omega_{n-1}\biggl(\frac{n-p}{p-1}\biggr)^{p-1}
\biggl(a^\frac{p-n}{p-1}-b^\frac{p-n}{p-1}\biggr)^{1-p},
\end{equation*}
where $\omega_{n-1}$ is the $(n-1)$-dimensional Lebesgue measure of
the unit sphere $S^{n-1}$ in $\mathbb{R}^n$ (see e.g. \cite{Geh71}).
Indeed, for $f(x)=\lambda x$, $\lambda>0$, $\lambda\in \mathbb R$,
we have $\mathcal M_p(f(\Gamma_A))=\lambda^{n-p}\mathcal
M_p(\Gamma_A)$.

\medskip
We also use especially another tool which is important in Potential
Theory and Mathematical Analysis.

Following in general \cite{MRV69}, a pair $E=(A,C)$, where
$A\subset\mathbb{R}^n$ is an open set and $C\subset A$ is a nonempty
compact, is called the \textit{condenser}. We say that the condenser
$E$ is the \textit{ring condenser}, if $\mathcal R=A\setminus C$ is
a ring domain. The condenser $E$ is bounded, if $A$ is bounded. We
also say that a condenser $E=(A,C)$ lies in a domain $G$ when
$A\subset G$. Obviously, for an open and continuous mapping
$f:G\to\mathbb{R}^n$ and for any condenser $E=(A,C)\subset G$, the
pair $(f(A),f(C))$ is a condenser in $f(G)$. In this case we shall
use the notation $f(E)=(f(A),f(C))$.

\medskip

Let $E=(A,C)$ be a condenser. Denote by $\mathcal C_0(A)$ the set of
all continuous functions $u:A\to\mathbb{R}^1$ with compact support
in $A$. Consider the set $\mathcal W_0(E)=\mathcal W_0(A,C)$ of all
nonnegative functions $u:A\to\mathbb{R}^1$ such that

\noindent 1) $u\in \mathcal C_0(A)$, 2) $u(x)\ge 1$ for $x\in C$ and
3) $u$ belongs ${\rm ACL}$. Put
\begin{equation*}
{\rm cap}_p\,E={\rm cap}_p\,(A,C)=\inf\limits_{u\in \mathcal
W_0(E)}\, \int\limits_{A}\,\vert\nabla u\vert^p\,dx, \quad p\ge 1,
\end{equation*}
where, as usual
\begin{equation*}
\vert\nabla u\vert={\left(\sum\limits_{i=1}^n\,{\left(\partial_i
u\right)}^2 \right)}^{1/2}.
\end{equation*}
This quantity is called \textit{$p$-capacity of condenser} $E$.

It was proven in \cite{Hes75} that for $p>1$
\begin{equation}\label{EMC}
{\rm cap}_p\,E=\mathcal M_p(\Delta(\partial A,\partial C; A\setminus
C)),\ \
 \end{equation}
where $\Delta(\partial A,\partial C; A\setminus C))$ denotes the set
of all continuous curves which join the boundaries $\partial A$ and
$\partial C$ in $A\setminus C$. For general properties of
$p$-capacities and their relation to the mapping theory, we refer
for instance to \cite{GR90} and \cite{Maz03}. In particular, for
$1\le p<n$,
\begin{equation}\label{maz}
{\rm cap}_p\,E\ge n{\Omega}^{\frac{p}{n}}_n
\left(\frac{n-p}{p-1}\right)^{p-1}\left[m C\right]^{\frac{n-p}{n}},
\end{equation}
where ${\Omega}_n$ denotes the volume of the unit ball in ${\Bbb
R}^n$, and $mC$ is the $n$-dimensional Lebesgue measure of $C$.

For $n-1<p\leq n$, there is the following lower estimate
\begin{equation}\label{krd}
\left({\rm cap}_{p}\,\,
E\right)^{n-1}\,\ge\,\gamma\,\frac{d(C)^{p}}{(mA)^{1-n+p}}\,\,,
\end{equation}
where $d(C)$ denotes the diameter of $C$, and $\gamma$ is a positive
constant depending only on $n$ and $p$ (see \cite {Krug86}).

\section{$Q$-homeomorphisms and their properties}
Let
$Q:\;G\to[1,\infty]$ be a measurable function. Due to \cite{MRSY04},
a homeomorphism $f:\;G\to \mathbb{R}^n$ is called a
\textit{$Q$-homeomorphism} if
\begin{equation}\label{Qh}
\mathcal M(f(\Gamma))\le\int\limits_G Q(x)\rho^n(x)\;dx
\end{equation}
for every family $\Gamma$ of curves in $G$ and for every admissible
function $\rho$ for $\Gamma$ (see also \cite{MRSY09}).

The origin of this notion relies on a natural generalization of
quasiconformality. Given a function $Q:\;G\to[1,\infty]$, we say
that a sense preserving homeomorphism $f:\;G\to \mathbb{R}^n$ is
{\it $Q(x)$-quasiconformal} if $f\in W^{1,n}_{\rm loc}(G)$ and
\begin{equation*}
\max\{H_I(x,f), H_O(x,f)\}\le Q(x) \quad\text{a.e.}
\end{equation*}

Any $Q(x)$-quasiconformal mapping $f:\;G\to \mathbb R^n$ is
differentiable a.e., satisfies the $(N)$-property, and $J(x,f)\ge 0$
a.e. If, in addition, $Q\in L^{n-1}_{\rm loc}$, then $f^{-1}\in
W_{\rm loc}^{1,n}(G^*)$ and is differentiable a.e.; $f$ has
$(N^{-1})$-property and $J(x,f)>0$ a.e. All this follows from
\cite{HK93, RR55, Re89} (cf. \cite{MRSY09}).

\medskip
Given a measurable function $Q:\;G\to[0,\infty]$, a homeomorphism
$f:\;G\to \mathbb R^n$ is called \textit{$Q$-homeomorphism with
respect to $\alpha$-module}, if
\begin{equation}\label{eq4.1}
\mathcal M_\alpha(f(\mathcal S_k))\le\int\limits_G
Q(x)\rho^\alpha(x)\;dx
\end{equation}
for every family of $k$-dimensional surfaces $\mathcal S_k$ in $G$
and for every admissible function $\rho$ for $\mathcal S_k$; $1\le
k\le n-1$, and such integer $k$ is fixed.

\medskip
It is the well-known fact that quasiconformal mappings preserve
their $n$-moduli up to an absolute factor, i.e.
\begin{equation}\label{eq3.1}
K^\frac{k-n}{n-1}\mathcal M(\mathcal S_k)\le \mathcal M(f(\mathcal
S_k))\le K^\frac{n-k}{n-1}\mathcal M(\mathcal S_k)
\end{equation}
(a quasiinvariance of $n$-module). For conformal mappings the
$n$-module becomes an invariant. Observe that the inequality
(\ref{Qh}) is a natural generalization of the right-hand side in
(\ref{eq3.1}) for the curve families. Note also that the integral in
(\ref{eq4.1}) can be interpreted as a weighted module (cf.
\cite{AC71, Cris11, Tam98}).

\medskip
We now restrict ourselves by the case $k=1$, which corresponds to
the curve families. For $Q$-homeomorphisms with $Q\in L_{\rm
loc}^1$, the differentiability a.e. and ACL-property were
established in \cite{Sal08}. It is also known that
$Q$-homeomorphisms satisfy the $(N^{-1})$-property (see
\cite{MRSY09}).

To establish the differential properties of $Q$-homeomorphisms with
respect to $\alpha$-moduli, we first consider some set functions.
Let $\Phi$ be a finite nonnegative function in $G$ defined for open
subsets $E$ of $G$ satisfying
\begin{equation*}
\sum_{k=1}^m\Phi(E_k)\le \Phi(E)
\end{equation*}
for any finite collection $\{E_k\}_{k=1}^m$ of nonintersecting open
sets $E_k\subset E$. Denote the class of all such set functions
$\Phi$ by $\mathcal F$.

The {\it upper} and {\it lower} derivatives of a set function
$\Phi\in \mathcal F$ at a point $x\in G$ are defined by
\begin{equation*}
\overline{\Phi'}(x)=\lim\limits_{h\to 0}\sup\limits_{d(Q)<h}
\frac{\Phi(Q)}{mQ},\quad\quad \underline{\Phi'}(x)=\lim\limits_{h\to
0}\inf\limits_{d(Q)<h} \frac{\Phi(Q)}{mQ},
\end{equation*}
where $Q$ ranges over all open cubes or open balls such that $x\in
Q\subset G$. Due to \cite {RR55}, these derivatives have the
following properties: $\overline{\Phi'}(x)$ and
$\underline{\Phi'}(x)$ are Borel's functions;
$\overline{\Phi'}(x)=\underline{\Phi'}(x)<\infty$ a.e. in $G$; and
for each open set $V\subset G$,
\begin{equation*}
\int\limits_V \overline{\Phi'}(x)\,dx\le\Phi(V).
\end{equation*}

\medskip
\begin{theo}\label{th4.1}
{\em Let $f:\;G\to G^*$ be an $Q$-homeomorphism with respect to
$\alpha$-module with $Q\in L^1_{\rm loc}(G)$ and $\alpha>n-1$. Then
$f$ is ACL-homeomorphism which is differentiable a.e. in $G$.}
\end{theo}

\medskip
For the proof of Theorem \ref{th4.1} we refer to \cite{Gol09} (cf.
\cite{Sal08} for $\alpha=n$).

\medskip
\begin{corol} Under the assumptions of Theorem \ref{th4.1},
any $Q$-homeomorphism with respect to $\alpha$-module belongs to
$W^{1,1}_{\rm loc}$.
\end{corol}

\medskip
The following theorem implies the upper estimates for the maximal
stretching and Jacobian of $f$.

\medskip
\begin{theo}
{\em Let $G$ and $G^*$ be two bounded domains in $\mathbb{R}^n$, $n\geq
2$, and $f:G\rightarrow G^*$ be a sense preserving $Q$-homeomorphism
with respect to $\alpha$-module, $n-1<\alpha<n$, so that $Q\in
L^1_{\rm loc}(G)$. Then
\begin{equation}\label{DF}
\|f'(x)\|\leq \lambda_{n,\alpha}\, Q^{\frac{1}{n-\alpha}}(x) \ \
\text{a.e.}
\end{equation}
and
\begin{equation}\label{Jacob}
J(x,f)\leq \gamma_{n,\alpha}\, Q^{\frac{n}{n-\alpha}}(x) \ \
\text{a.e.},
\end{equation}
where $\lambda_{n,\alpha}$ and $\gamma_{n,\alpha}$ are  a positive
constants depending only on $n$ and $\alpha$.}
\end{theo}

\begin{proof}
First consider the set function $\Phi(B)=mf(B)$ defined over the
algebra of all the Borel sets $B$ in $G$. By \cite{MRV69, RR55},
\begin{equation*}
\varphi(x)=\limsup\limits_{\varepsilon\rightarrow
0}\frac{\Phi(B(x,\varepsilon))}{\Omega_{n}\varepsilon^{n}}<\infty
\quad\text{for a.e}\quad x\in G.
\end{equation*}

Let $A=A(x_0,\varepsilon_1, \varepsilon_2)=\{x:\
\varepsilon_1<|x-x_0|< \varepsilon_2\}$ be a spherical ring centered
at $x_0\in G,$ with radii $\varepsilon_1$ and $\varepsilon_2$,
$0<\varepsilon_1<\varepsilon_2$, such that $A(x_0,\varepsilon_1,
\varepsilon_2)\subset G$. Then
$\left(f(B(x_0,\varepsilon_2)),f(\overline{B}(x_0,\varepsilon_1))\right)$
is a ring condenser in  $G^*$ and in accordance with (\ref{EMC}),
\begin{equation*}
{\rm cap}_{\alpha}\,(f(B(x_0,\varepsilon_2)),f(\overline{B}(x_0
,\varepsilon_1))=\mathcal M_{\alpha}(\triangle(\partial
f(B(x_0,\varepsilon_2)),\partial f(B(x_0,\varepsilon_1));f(A)).
\end{equation*}
Since $f$ is homeomorphic
\begin{equation*}
\triangle\left(\partial f(B(x_0,\varepsilon_2)),\partial
f(B(x_0,\varepsilon_1));f(A)\right)=f\left(\triangle\left(\partial
B(x_0,\varepsilon_2) ,\partial B(x_0,\varepsilon_1);A\right)\right).
\end{equation*}

Pick the admissible function
\begin{equation*}
\eta(t)\,=\,\left
\{\begin{array}{rr} \frac{1}{\varepsilon_2-\varepsilon_1}, &  \ t\in (\varepsilon_1,\varepsilon_2) \\
0, & \ t\in \Bbb{R}\setminus (\varepsilon_1,\varepsilon_2).
\end{array}\right.
\end{equation*}
Since $f$ is a $Q$-homeomorphism with respect to $\alpha$-module,
\begin{equation}\label{eq100}
{\rm
cap}_{\alpha}\,(f(B(x_0,\varepsilon_2)),f(\overline{B}(x_0,\varepsilon_1))
\leq
\frac{1}{(\varepsilon_2-\varepsilon_1)^{\alpha}}\int\limits_{A(x_0,\varepsilon_1,
\varepsilon_2)} Q(x)\, dx.
\end{equation}
Choosing $\varepsilon_1=2\varepsilon$ and
$\varepsilon_2=4\varepsilon$, we get
\begin{equation}\label{eq101}
{\rm
cap_{\alpha}}\,(f(B(x_0,4\varepsilon)),f(\overline{B}(x_0,2\varepsilon))\le\,
\frac{1}{(2\varepsilon)^{\alpha}}\int\limits_{B(x_0,4\varepsilon)}Q(x)\,dx.
\end{equation}
On the other hand, the inequality (\ref{maz}) implies
\begin{equation}\label{eq102}
{\rm
cap_{\alpha}}\,(fB(x_0,4\varepsilon),f(\overline{B}(x_0,2\varepsilon))
\ge
C_{n,\alpha}\left[mf(B(x_0,2\varepsilon))\right]^{\frac{n-\alpha}{n}},
\end{equation}
where the constant $C_{n,\alpha}$ depends only on the dimension $n$
and $\alpha.$

Combining (\ref{eq101}) and  (\ref{eq102}) yields
\begin{equation}\label{eq4.2}
\frac{mf(B(x_0,2\varepsilon))}{mB(x_0,2\varepsilon)}\le
c_{n,\alpha}\,\left[ \dashint_{B(x_0,4\varepsilon)}\, Q(x)\, dx
\right]^{\frac{n}{n-\alpha}},
\end{equation}
where $c_{n,\alpha}$ also depends on $n$ and $\alpha$. As
$\varepsilon\to 0$, the estimate (\ref{Jacob}) follows.

\medskip
Now choosing $\varepsilon_1=\varepsilon$ and
$\varepsilon_2=2\varepsilon$ in (\ref{eq100}), we obtain
\begin{equation*}
{\rm cap_{\alpha}}\
(f(B(x_0,2\varepsilon)),f(\overline{B}(x_0,\varepsilon))\le\,
\frac{1}{\varepsilon^{\alpha}}\int\limits_{B(x_0,2\varepsilon)}Q(x)\,dx,
\end{equation*}
and since by (\ref{krd}),
\begin{equation*}
{\rm cap_{\alpha}}\
(f(B(x_0,2\varepsilon)),f(\overline{B}(x_0,\varepsilon)) \ge
C_{n,\alpha}\frac{[d(f(B(x_0,\varepsilon)))]^\frac{\alpha}{n-1}}{[mf(B(x_0,2\varepsilon))]^\frac{1-n+\alpha}{n-1}},
\end{equation*}
we have
\begin{equation}\label{eq92}
\frac{d(fB(x_0,\varepsilon))}{\varepsilon}\le\gamma_{n,\alpha}
\left(\frac{mf(B(x_0,2\varepsilon))}{mB(x_0,2\varepsilon)}
\right)^{\frac{1-n+\alpha}{\alpha}}\left(
\dashint_{B(x_0,2\varepsilon)} Q(x)\,
dx\right)^{\frac{n-1}{\alpha}}\,,
\end{equation} where $\gamma_{n,\alpha}$ is a positive constant depending on  $n$ and $\alpha.$

The inequalities (\ref{eq92}) and (\ref{eq4.2}) yield
\begin{equation*}
\frac{d(f(B(x_0,\varepsilon)))}{\varepsilon}\leq\lambda_{n,\alpha}\left(\dashint_{B(x_0,4\varepsilon)}
Q(x)\, dx\right)^{\frac{n(1-n+\alpha)}{\alpha(n-\alpha)}}\left[
\dashint_{B(x_0,2\varepsilon)}\, Q(x)\, dx
\right]^{\frac{n-1}{\alpha}}\,.
\end{equation*}
Finally, letting $\varepsilon\to 0$, one derives
\begin{equation*}
\limsup\limits_{x\to
x_0}\frac{|f(x)-f(x_0)|}{|x-x_0|}\leq\limsup\limits_{\varepsilon\to
0}\frac{d(f(B(x_0,\varepsilon)))}{\varepsilon}\leq
\lambda_{n,\alpha}\, Q^{\frac{1}{n-\alpha}}(x_{0}) ,
\end{equation*}
where  $\lambda_{n,\alpha}$ is a positive constant which depends
only on $n$ and $\alpha$. Thus (\ref{DF}) follows.

\end{proof}

\begin{corol}
Let  $G$ and $G^*$ be two domains in $\mathbb{R}^n$, $n\ge 2$, and
let $f:G\to G^*$ be a $Q$-homeomorphism with respect to
$\alpha$-module, $n-1<\alpha<n$. Assume that $Q(x)\in
L^{\frac{s}{n-\alpha}}_{\rm loc}(G)$ with $s>n-\alpha$. Then $f\in
W^{1,s}_{\rm loc} $.
\end{corol}

\medskip
Indeed for any compact set $V\subset G$,
\begin{equation*}
\int\limits_{V}\|f'(x)\|^s\ dx\leq
\lambda^s_{n,\alpha}\int\limits_{V}Q^{\frac{s}{n-\alpha}}(x)\
dx\,<\infty\,.
\end{equation*}

\medskip
As well known, every $W^{1,n}_{\rm loc}$-homeomorphism possesses the
$(N)$-property; thus we have

\medskip
\begin{corol}
If $Q\in L^\frac{n}{n-\alpha}_{\rm loc}(G)$ then $f$ satisfies
$(N)$-property.
\end{corol}

\medskip
\begin{corol}
Let $f:G\to G^*$ be a $Q$-homeomorphism with respect to
$\alpha$-module such that $Q(x)\in L^{\frac{n}{n-\alpha}}_{\rm
loc}(G)$, $n-1<\alpha<n$. Then
\begin{equation*}
mf(E)\leq  \gamma_{n,\alpha} \int\limits_{E}
Q^{\frac{n}{n-\alpha}}(x)\, dx.
\end{equation*}
\end{corol}

\begin{proof}
Since $Q(x)\in L^{\frac{n}{n-\alpha}}_{\rm loc}(G)$, $f$ satisfies
Lusin's $(N)$-property and
\begin{equation*}
mf(E)=\int\limits_{E} \, J(x,f)\, dx\leq
\gamma_{n,\alpha}\int\limits_{E} Q^{\frac{n}{n-\alpha}}(x)\, dx\,.
\end{equation*}
\end{proof}

The following result shows that the dilatation $H_{I,\alpha}$ is
dominated a.e. in $G$ by the upper derivative of the set function
\begin{equation*}
\Psi(V)=\int\limits_V Q(x)\,dx,
\end{equation*}
where $V$ is an open subset of $G$.

\medskip
\begin{theo}\label{th4.2}
{\em Let $f:\;G\to G^*$ be a $Q$-homeomorphism with respect to
$\alpha$-module, $\alpha>n-1$, with $Q\in L^1_{\rm loc}(G)$ such that
$J(x,f)\ne 0$ a.e. in $G$. Then
\begin{equation*}
H_{I,\alpha}(x,f)\le \overline\Psi'(x)\quad\text{a.e. in}\quad G.
\end{equation*}}
\end{theo}

\begin{proof}
Let $a\in G$ be an arbitrary point where $f$ is differentiable at
$a$, with $J(a,f)\ne 0$ and $\overline\Psi'(a)\ne 0$. The image of the unit ball under the linear
mapping $f^\prime(a)$ is an ellipsoid $E(f)$ with semi-axes
$\l_1,\l_2,\ldots,\l_n$ ordered by $\l_1\ge\l_2\ge\ldots\ge\l_n>0$.
Preceding $f$, if necessary, by a rotation and a reflection, one
reduces to the case $f(a)=a=0$ and $|f'(0)e_i|=\l_i$,
$i=1,\ldots,n$; here $e_\nu$ denotes the $\nu$th unit basis vector.

For every $t>0$, let $\mathcal R$ be the ring domain obtained from
$n$-dimensional interval
\begin{equation*}
I_n=\{x: |x_i|<r(t\l_i+1), i=1,\ldots,n-1, |x_n|<rt\l_n\},
\end{equation*}
by deleting the points of $(n-1)$-dimensional interval
\begin{equation*}
\Pi_{n-1}(0,r)=\{x:|x_i|\le r, i=1,\ldots,n-1, x_n=0\}.
\end{equation*}
We choose $r>0$ so that $\overline {\mathcal R}\subset G$ and will
show that
\begin{equation}\label{eq5.1}
\frac{\l_1\cdot\ldots\cdot\l_n}{\l_n^\alpha} \le \overline\Psi'(0).
\end{equation}
Indeed, the inequality
\begin{equation*}
\mathcal M_\alpha(f(\Gamma_\mathcal R))\le\int\limits_\mathcal R
Q(x)\rho^\alpha(x)\;dx,
\end{equation*}
together with the following estimate (see, e.g. \cite {Krug86})
\begin{equation*}
\mathcal M_\alpha(\Gamma_\mathcal R)\ge
\frac{\bigl(\inf\limits_{\Sigma} m_{n-1}
\S\bigr)^\alpha}{\bigl(m\mathcal R\bigr)^{\alpha-1}},
\end{equation*}
gives
\begin{equation*}
\frac{\biggl(\inf\limits_{\Sigma^*} m_{n-1}
\Sigma^*\biggr)^\alpha}{\bigl(mf(\mathcal R)\bigr)^{\alpha-1}} \le
\frac{1}{\bigl(\dist(F_0,F_1)\bigr)^\alpha} \int\limits_\mathcal R Q(x)\,dx,
\end{equation*}
where the infimum is taken over all surfaces $\Sigma^*$ which
separate $f(C_0)$ and $f(C_1)$ in $f(\mathcal R)$.

Fix $0<\varepsilon<\l_n$ and choose $r>0$ so small that
\begin{equation*}
mf(I_n)\le(J(0,f)+\varepsilon)mI_n \quad \text {and} \quad
|f(x)-f'(0)x|<\varepsilon r.
\end{equation*}
Since
\begin{equation*}
m\mathcal R=mI_n=2^n r^n
t\l_n(t\l_1+1)\cdot\ldots\cdot(t\l_{n-1}+1),
\end{equation*}
and
\begin{equation*}
\inf\limits_{\Sigma^*} m_{n-1} \Sigma^*\ge 2m_{n-1}\tilde
\Pi_{n-1}(0,r)=2^n
r^{n-1}(\l_1-\varepsilon)\cdot\ldots\cdot(\l_{n-1}-\varepsilon),
\end{equation*}
we have
\begin{equation*}
\frac{\bigl[2^nr^{n-1}(\lambda_1-\varepsilon)\cdot\ldots\cdot
(\lambda_{n-1}-\varepsilon)\bigr]^\alpha}
{\bigl[(J(0,f)+\varepsilon)mA\bigr]^{\alpha-1}}\le
\frac{1}{\bigl(rt\lambda_n\bigr)^\alpha} \int\limits_\mathcal R Q(x)\,dx,
\end{equation*}
where $\tilde \Pi_{n-1}(0,r)=\{y:|y_i|\le r\l_i-\varepsilon r,
i=1,\ldots,n-1, y_n=0\}$. Letting $t\to 0$ and thereafter $r\to 0$,
we get
\begin{equation*}
\frac {\bigl[(\l_1-\varepsilon)\cdot\ldots
\cdot(\l_{n-1}-\varepsilon)\bigr]^\alpha}
{\bigl[J(0,f)+\varepsilon\bigr]^{\alpha-1}} \le\overline\Psi'(0),
\end{equation*}
which implies (\ref{eq5.1}) as $\varepsilon\to 0$. Hence
$H_{I,\alpha}(x,f)\le \overline\Psi'(x)$ for almost all $x\in
\mathcal R$.
\end{proof}

\begin{rem} If one omits the restriction $\alpha>n-1$ (i.e. for $1\le\alpha\le n-1$),
Theorem \ref{th4.2} can be proved assuming additionally that the
$Q$-homeomorphisms $f$ with respect to $\alpha$-module are
differentiable a.e.
\end{rem}

\medskip
\begin{rem}
The inequality in Theorem \ref{th4.2} can be replaced by $H_{I,\alpha}(x,f)\le Q(x)$ a.e.
\end{rem}

\section{Ring $Q$-homeomorphisms and their properties}
Recall
some necessary notions. Let $E\,,F\,\subseteq\,\mathbb R^n$ be
arbitrary domains. Denote by $\Delta(E,F,G)$ the family of all
curves $\gamma:[a,b]\,\rightarrow\,\mathbb R^n$, which join $E$ and
$F$ in $G$, i.e. $\gamma(a)\,\in\,E\,,\gamma(b) \,\in\,F$ and
$\gamma(t)\,\in\,G\,$ for $a\,<\,t\,<\,b\,.$ Set $d_0\,=\,{\rm
dist}\,(x_0\,,\partial G)$ and let
$Q:\,G\rightarrow\,[0\,,\infty]\,$ be a Lebesgue measurable
function. Denote
\begin{equation*}
A(x_0,r_1,r_2) = \{ x\,\in\,{\mathbb R}^n : r_1<|x-x_0|<r_2\}\ ,
\end{equation*}
and
\begin{equation}\label{2.11}
S_i\,=\,S(x_0,r_i) = \{ x\,\in\,{\Bbb R^n} : |x-x_0|=r_i\}\,\,,\ \ \
i=1,2.
\end{equation}

\medskip
We say that a homeomorphism $f:G\to \mathbb R^n$  is \textit{the
ring $Q$-homeomorphism with respect to $p$-module at the point
$x_0\,\in\,G,$} ($1<p\leq n$) if the inequality
\begin{equation}\label{2.12}
\mathcal M_{p}\,\left(\Delta\left(f(S_1),f(S_2),
f(G)\right)\right)\,\ \leq \int\limits_{A} Q(x)\cdot
\eta^{p}(|x-x_0|)\ d\,x
\end{equation}
is fulfilled for any ring $A=A( x_0,r_1,r_2),$\,\, $0<r_1<r_2< d_0$
and for every measurable function $\eta : (r_1,r_2)\to [0,\infty
]\,,$ satisfying
\begin{equation}\label{2.13}
\int\limits_{r_1}^{r_2}\eta(r)\ dr\ \geq\ 1\,.
\end{equation}
The homeomorphism  $f:G\to \mathbb R^n$ is the \textit{ring
$Q$-homeomorphism with respect to $p$-module in the domain $G$}, if
inequality $(\ref{2.12})$ holds for all points $x_0\,\in\,G\,.$ The
properties of the ring $Q$-homeomorphisms for $p=n$ are studied in
\cite{SS10}.

The ring $Q$-homeomorphisms are defined in fact locally and contain
as a proper subclass of $Q$-homeomorphisms (see \cite{MRSY09}). A
necessary and sufficient condition for homeomorphisms to be ring
$Q$-homeomorphisms with respect to $p$-module at a point given in
\cite{Sal11}, asserts:

\medskip
\begin{propo}\label{th5.1}
{\em Let $G$ be a bounded domain in $\mathbb R^n$, $n\ge 2$ and let
$Q:\,G\,\rightarrow\,[0,\,\infty]\,$ belong to $L^1_{\rm loc}(G)$. A
homeomorphism $f:\,G\to \mathbb R^n$ is a ring $Q$-homeomorphism
with respect to $p$-module at $x_0\in G$ if and only if for any
$0<r_1<r_2< d_0= dist\,(x_0,\partial G)$,
\begin{equation*}
\mathcal M_p\left(\Delta\left(f(S_1),\,f(S_2),\,
f(G)\right)\right)\,\le\,\frac{\omega_{n-1}}{I^{p-1}}\,\,,
\end{equation*}
where $S_1$ and $S_2$ are the spheres defined in (\ref{2.11})
\begin{equation*}
I\ =\ I(x_0,r_1,r_2)\ =\ \int\limits_{r_1}^{r_2}\
\frac{dr}{r^{\frac{n-1}{p-1}} q_{x_0}^{\frac{1}{p-1}}(r)}\,,
\end{equation*}
and $q_{x_0}(r)\,$ is the mean value of $Q$ over $|x-x_0|\,=\,r\,.$
Note that the infimum in the right-hand side of $(\ref{2.12})$ over
all admissible $\eta$ satisfying (\ref{2.13}) is attained only for
the function
\begin{equation*}
\eta_0(r)=\frac{1}{Ir^{\frac{n-1}{p-1}}
q_{x_0}^{\frac{1}{p-1}}(r)}\,.
\end{equation*}}
\end{propo}

\section{Lower $Q$-homeomorphisms and their module bounds}
Let
$G$ and $G^*$ be two bounded domains in $\mathbb R^n$, $n\ge 2$ and
$x_0\in G$. Given a Lebesgue measurable function $Q:G\to
[0,\infty]$, a homeomorphism $f:G\to G^*$ is called the
\textit{lower $Q$-homeomorphism with respect to $p$-module at} $x_0$
if
\begin{equation}\label{eq6.1}
\mathcal M_p(f(\Sigma_\varepsilon))\ge \inf\limits_{\rho\in{\rm
expadm}\;\Sigma_\varepsilon}
\int\limits_{D_\varepsilon(x_0)}\frac{\rho^p(x)}{Q(x)}dx,
\end{equation}
where
\begin{equation*}
G_\varepsilon(x_0)=G\cap \{x\in\mathbb R^n:
\varepsilon<|x-x_0|<\varepsilon_0\},\quad
0<\varepsilon<\varepsilon_0,\quad 0<\varepsilon_0<\sup\limits_{x\in
G}|x-x_0|,
\end{equation*}
and $\Sigma_\varepsilon$ denotes the family of all pieces of spheres
centered at $x_0$ of radii $r$, $\varepsilon<r<\varepsilon_0$,
located in $G$.

In this section, we provide a necessary and sufficient condition for
homeomorphisms to be lower $Q$-homeomorphisms with respect to
$p$-module. The case $p=n$ is considered in \cite{MRSY09}.

\medskip
\begin{theo}\label{th6.1}
{\em Let $G$ be a domain in $\mathbb{R}^n$, $n\geq2$,
$x_0\in\overline{G}$, and let $Q:\,G\to[0,\infty]$ be a measurable
function. A homeomorphism $f:\,G\to\mathbb{R}^n$ is a lower
$Q$-homeomorphism at $x_0$ with respect to $p$-module for $p>n-1$ if
and only if the following inequality
\begin{equation*}
\mathcal
M_p(f(\Sigma_{\varepsilon}))\geq\int\limits_{\varepsilon}^{\varepsilon_0}
\frac{dr}{||\,Q||\,_{s}(r)}\quad\forall\
\varepsilon\in(0,\varepsilon_0)\,,\quad\varepsilon_0\in(0,d_0)\,,
\end{equation*}
holds, where $s=\frac{n-1}{p-n+1}$,
\begin{equation*}
d_0=\sup\limits_{x\in D}\,|x-x_0|\,,
\end{equation*}
$\Sigma_{\varepsilon}$ is the family of all intersections of the
spheres $S(x_0,r)=\{x\in\mathbb{R}^n:|\,x-x_0|=r\}$,
$r\in(\varepsilon,\varepsilon_0)$ with $G$ and
\begin{equation*}
||Q||_{s}(r)=\left(\int\limits_{G(x_0,r)}
Q^{s}(x)\,d\sigma_{n-1}\right)^\frac{1}{s}\;
\end{equation*}
here $G(x_0,r)=\{x\in G:|\,x-x_0|=r\}=G\cap S(x_0,r)$. The infimum
in (\ref{eq6.1}) is attained only on the functions
\begin{equation*}
\varrho_0(x)=\frac{Q(x)}{||Q||_s\left(r\right)}\,.
\end{equation*}}
\end{theo}

\begin{proof}
For any $\varrho\in {\rm extadm}\,\Sigma_{\varepsilon}$ the function
\begin{equation*}
A_{\varrho}(r):=\int\limits_{G(x_0,r)}\varrho^{n-1}(x)\,d\sigma_{n-1}\neq0\qquad
\end{equation*}
is measurable on $(\varepsilon,\varepsilon_0)$. It is admissible if
$A_{\varrho}(r)\ge 1$. Assuming that $A_{\varrho}(r)\equiv 1$, we
obtain
\begin{equation*}
\inf\limits_{\varrho\in {\rm extadm}\,\Sigma_{\varepsilon}}
\int\limits_{G_{\varepsilon}(x_0)}\frac{\varrho^p(x)}{Q(x)}\,dx=
\int\limits_{\varepsilon}^{\varepsilon_0}\left(\inf\limits_{\psi\in
I(r)}\int\limits_{G(x_0,r)}\frac{\psi^q
(x)}{Q(x)}\,d\sigma_{n-1}\right)dr\,,
\end{equation*}
where $q=p/(n-1)>1$ and $I(r)$ denotes the set of all measurable
functions $\psi$ on the surface $G(x_0,r)$ satisfying
\begin{equation*}
\int\limits_{G(x_0,r)}\psi(x)\,d\sigma_{n-1}=1\,.
\end{equation*}

Thus Theorem \ref{th6.1} follows from Lemma \ref{lem} taking there
$X=G(x_0,r)$ and $\mu$ to be the $(n-1)$-dimensional area on
$G(x_0,r)$, and $\varphi(x)=1/Q(x)$ on $G(x_0,r)$, and
$q=p/(n-1)>1$. This completes the proof.
\end{proof}

\section{Connection between the ring and lower
$Q$-homeomorphisms}
In this sections we establish the relationship
between the ring and lower $Q$-homeomorphisms with respect to
$p$-module.

\medskip
\begin{theo}\label{th7.1}
{\em Every lower $Q$-homeomorphism with respect to $p$-module $f:G\rightarrow G^*$ at $x_0\in G,$ with $p>n-1$ and
  $Q\in L_{loc}^{\frac{n-1}{p-n+1}}(G)$ is a ring
$\widetilde{Q}$-homeomorphism with respect to $\alpha$-module at
$x_0$ with $\widetilde{Q}=Q^{\frac{n-1}{p-n+1}}$ and
$\alpha=\frac{p}{p-n+1}$.}
\end{theo}

\begin{proof} Let
$0<r_1<r_2<d(x_0, \partial G)$ and $S_i=S(x_0, r_i)$, $i=1,2$, be
from (\ref{2.11}). Then taking into account the relation (\ref{ziem}), we obtain
\begin{equation}\label{eqOS6.1aa}
\mathcal M_\alpha\left(f\left(\Delta(S_{1}, S_{2}, G)\right)\right)\
\le \frac{1}{\mathcal
M_p^{\frac{n-1}{p-n+1}}(f\left(\Sigma\right))}\,,
\end{equation}
where
$f\left(\Sigma\right)\subset\Sigma\left(f(S_1),f(S_2),f(G)\right)$
and $\Sigma$ denotes the family of all spheres centered at $x_0,$
located between $S_1$ and $S_2,$ while
$\Sigma\left(f(S_1),f(S_2),f(G)\right)$ consists of all
$(n-1)$-dimensional surfaces in $f(G)$, separating $f(S_1)$ and
$f(S_2)$ (cf. \cite{MRSY09}) Now directly by (\ref{eqOS6.1aa}) and
Theorem \ref{th6.1},
\begin{equation}\label{eqOS6.1b}
\mathcal M_\alpha\left(f\left(\Delta(S_{1}, S_{2}, G)\right)\right)\
\le \left(\int\limits_{r_1}^{r_2}
\frac{dr}{\Vert\,Q\Vert\,_{\frac{n-1}{p-n+1}}(r)}\right)^{\frac{1-n}{p-n+1}}=
\frac{\omega_{n-1}}{\widetilde{I}^\frac{n-1}{p-n+1}}\,.
\end{equation}
where $\widetilde{I}\ =\ \widetilde{I}(x_0,r_1,r_2)\ =\
\int\limits_{r_1}^{r_2}\ \frac{dr}{r^{\frac{n-1}{\alpha-1}}
\widetilde{q}_{x_0}^{\frac{1}{\alpha-1}}(r)}\,,$ and
$\widetilde{q}_{x_0}(r)\,$ denotes the mean value of the function
$\widetilde{Q}$ over $|x-x_0|\,=\,r\,.$

Now the assertion of the theorem follows from (\ref{eqOS6.1b}) and
Proposition \ref{th5.1}.
\end{proof}

\medskip
{\small \leftline{\textbf{Anatoly Golberg}} \em{
\leftline{Department of Applied Mathematics,} \leftline{Holon
Institute of Technology,} \leftline{52 Golomb St., P.O.B. 305,}
\leftline{Holon 58102, ISRAEL} \leftline{Fax: +972-3-5026615}
\leftline{e-mail: golberga@hit.ac.il}}}

\medskip

{\small \leftline{\textbf{Ruslan Salimov}}\em \leftline{Institute of
Applied Mathematics and Mechanics,}\leftline {National Academy of
Sciences of Ukraine,} \leftline{74 Roze Luxemburg St.,}
\leftline{Donetsk 83114, UKRAINE}\leftline{e-mail:
salimov07@rambler.ru}}}

\end{document}